\documentstyle[11pt]{article}
\newcommand{\qed}{\hfill\rule{4pt}{8pt}\par\vspace{\baselineskip}}

\newtheorem{de}{Definition}[section]

\newtheorem{pr}[de]{Proposition}
\newtheorem{co}[de]{Corollary}
\newtheorem{re}[de]{Remark}

\newtheorem{te}[de]{Theorem}

\def\Box{\mbox{$\sqcap\!\!\!\!\sqcup$}}

\def\ot{\otimes}

\def\eps{\varepsilon}

\def\bea{\begin{eqnarray*}}
\def\eea{\end{eqnarray*}}

\begin{document}
\title{Frobenius Structural Matrix Algebras}
\author{
S. D\u{a}sc\u{a}lescu$^1$, M. C. Iovanov$^{1,2}$ and S. Predu\c{t}$^1$ \\[2mm]
$^1$University of Bucharest, Facultatea de Matematic\u{a},\\ Str.
Academiei 14, Bucharest 1, RO-010014, Romania,\\
$^2$ University of Iowa,
14 MacLean Hall, Iowa City, Iowa 52242, USA\\
e-mail: sdascal@fmi.unibuc.ro, miodrag-iovanov@uiowa.edu,
sina.predut@gmail.com}

\date{}
\maketitle

\begin{abstract}
We discuss when the incidence coalgebra of a locally finite
preordered set is right co-Frobenius. As a consequence, we obtain
that a structural matrix algebra over a field $k$ is Frobenius if
and only if it consists, up to a permutation of rows and columns,
of diagonal blocks which are full matrix
algebras over $k$.
\end{abstract}

\section{Introduction and preliminaries}

Let $k$ be a field. A structural matrix algebra is an algebra
consisting of all matrices in some $M_n(k)$ which have zeros on
certain prescribed positions. This kind of algebras are present in
ring theory by the large number of examples and counterexamples
using them. Also, certain structural matrix algebras are important
in the study of numerical invariants of PI algebras. Structural
matrix algebras have been studied in for example \cite{vw},
\cite{c}, \cite{abw}.

On the other hand, Frobenius algebras have been a theme of
research with applications in representation theory, topology,
quantum Yang-Baxter equation, conformal field theory, Hopf algebra
theory; see \cite{kadison}, \cite{lam}, \cite{sy}. It is a natural
question to ask when is a structural matrix algebra Frobenius. We
will answer this question in the present note. In fact we will
investigate a more general situation. We consider a locally finite
preordered set $(X,\leq )$, i.e. the relation $\leq$ is reflexive
and transitive, and the interval $[x,y]=\{z|\; x\leq z\leq y\}$ is
finite for any $x\leq y$. We consider the vector space $C=IC(X)$
with basis a set $\{ e_{x,y}|x,y\in X, x\leq y\}$. This space has
a coalgebra structure with comultiplication $\Delta$ and counit
$\eps$ defined on basis elements by
$$\Delta(e_{x,y})=\sum_{x\leq z\leq y}e_{x,z}\ot e_{z,y}$$
$$\eps (e_{x,y})=\delta_{x,y}$$
for any $x,y\in X$ with $x\leq y$; here  $\delta_{x,y}$  denotes
Kronecker's delta. Incidence coalgebras have been of a great
interest for combinatorics, see for instance \cite{jr}, \cite{so}.
We will discuss when is this incidence coalgebra right
co-Frobenius, i.e. when $C$ embeds in $C^*$ as a right
$C^*$-module. Co-Frobenius coalgebras have a rich representation
theory, also their study is interesting since a Hopf algebra has
non-zero integrals if and only if its underlying coalgebra is
right co-Frobenius. A similar problem was posed in \cite{dnv} in
the particular case of incidence coalgebras of locally finite
ordered sets. In that paper, an equivalent characterization of
right co-Frobenius coalgebras is considered,
 more precisely the existence of a right non-singular
$C^*$-balanced bilinear form on $C$, and the method is to compute
all such bilinear forms. We propose here a completely different
approach, by using corepresentation (i.e. comodule) theory of $C$.
Since a finite dimensional algebra is Frobenius if and only if its
dual coalgebra is right (or left) co-Frobenius, and the structural
matrix algebras are precisely the incidence algebras of finite
preordered sets (so they are isomorphic to the dual algebras of
incidence coalgebras of finite preordered sets), we obtain as an
application the structure of all Frobenius structural matrix
algebras.  We refer to \cite{dnr} for basic facts about coalgebras
and comodules. We will freely regard a right $C$-comodule as a
left $C^*$-module in the usual way.

\section{The coradical filtration}

We start with a general remark on coalgebras having a special kind
of basis. We will then apply it to incidence coalgebras. Let $C$
be an arbitrary coalgebra with comultiplication $\Delta$. The
coradical $C_0$ of $C$ is the sum of all simple subcoalgebras of
$C$, and this coincides with the socle of $C$ as a right (or as a
left) $C$-comodule. If $U$ and $V$ are subspaces of  $C$, the
wedge $U\wedge V$ is  defined by $U\wedge V=\Delta^{-1}(U\ot
C+C\ot V)$. The coradical filtration $C_0\subseteq C_1\subseteq
\ldots \subseteq C_n\subseteq \ldots$ of $C$ is defined by
$C_n=C_0\wedge C_{n-1}$ for any $n\geq 1$.

\begin{de} \label{defcoradbasis}
(i) We say that $\mathcal{B}$ is a coradical basis of $C$ if
$\mathcal{B}$
is a $k$-basis of $C$ and ${\mathcal{B}}\cap C_n$ is a basis of $C_n$ for all $n$.\\
(ii) A right (or left) coideal $M$ of $C$ will be called
$\mathcal{B}$-supported if ${\mathcal{B}}\cap M$ is a $k$-basis of
$M$.
\end{de}

We say that a non-zero (co)module is local if it has a unique
maximal sub-(co)module. The following result is obvious.

\begin{pr}\label{common}
Let $C$ be a coalgebra with a coradical basis $\mathcal{B}$ and
let $M$ be a finite dimensional local right coideal of $C$ which
is $\mathcal{B}$-supported. Then $M$ is generated as a right
$C$-comodule by any element in the nonempty set
$({\mathcal{B}}\cap M)\setminus Jac(M)$ (here, $Jac(M)$ denotes
the Jacobson radical of $M$).
\end{pr}

Let now $C$ be the incidence coalgebra $IC(X)$ of a locally finite
preordered set $X$.

Let $\sim$ be the equivalence relation on $X$ defined by $x\sim y$
if and only if $x\leq y$ and $y\leq x$. We denote by $({\cal
C}_i)_{i\in I}$ the set of equivalence classes with respect to
$\sim$, which are all finite since $\leq$ is locally finite. Then
for any $i\in I$, the space $D_i=\sum _{x,y\in {\cal
C}_i}ke_{x,y}$ is a subcoalgebra of $C$ isomorphic to the matrix
coalgebra $M^c(n_i,k)$, the dual coalgebra of the matrix algebra
$M_{n_i}(k)$. In particular $D_i$ is a simple coalgebra. If $i\in
I$ and $x\in {\cal C}_i$, define $S_x=\sum_{y\in {\cal
C}_i}ke_{x,y}$. Then $S_x$ is a simple right $C$-comodule and
$D_i=\sum_{x\in {\cal C}_i}S_x$.

For any $x\in X$ we also consider the space $E_x=<e_{x,y}\; |\;
y\in X, x\leq y >$, which is a right $C$-subcomodule of $C$.
Obviously $S_x\subseteq E_x$ and $C=\oplus_{x\in X}E_x$, thus any
$E_x$ is an injective right $C$-comodule.

If $x\leq y$, we say that the length of the interval $[x,y]$ is 0
if $x\sim y$, and the length of $[x,y]$ is $n>0$ if $n$ is the
greatest possible such that there exists a sequence
$x=x_0<x_1<\ldots <x_n=y$ in $X$.

\begin{pr}
(i) The coradical $C_0$ of $C$ is $\sum_{i\in I}D_i$.\\
(ii) $C_n$, the $(n+1)$th term of the coradical filtration of $C$,
is the subspace spanned by all $e_{x,y}$ such that $[x,y]$ has
length at most $n$.
\end{pr}
{\bf Proof:} Denote $D=\sum_{i\in I}D_i$. Since $D$ is a sum of
simple subcoalgebras, we see that $D\subseteq C_0$.

 It is straightforward to show by induction on $n$
that $\wedge^{n+1}D$ is the subspace spanned by all $e_{x,y}$ such
that $[x,y]$ has length at most $n$. It follows that
$\bigcup_n(\wedge^{n+1}D)=C$, where $\wedge^{n+1}D$ is the space
obtained by wedging $n$ times $D$ with itself. Then $C_0\subseteq
D$ by \cite[Exercise 3.1.12]{dnr}. We conclude that $C_0=D$. This
clearly implies (ii). \qed

\begin{co}
Any simple right $C$-comodule is isomorphic to $S_x$ for some
$x\in X$.
\end{co}

\begin{co}
The set of all $e_{x,y}$'s with $x\leq y$ is a coradical basis of
$C$.\end{co}

For any $x\leq y$ in $X$ we consider the element $p_{x,y}\in C^*$
such that $p_{x,y}(e_{u,v})=\delta_{x,u}\delta_{y,v}$ for any
$u,v\in X$ with $u\leq v$. If we denote by $\rightarrow$ the left
$C^*$-action on $C$ induced by the right $C$-comodule structure,
it is easy to check that

$$p_{x,y}\rightarrow e_{u,v}=\left\{
\begin{array}{l}
e_{u,x}, \mbox{ if } y=v \mbox{ and }u\leq x\\
0, \mbox{ \rm otherwise}
\end{array}\right.$$

\begin{pr}
For any $x\in X$, $S_x$ is an essential subcomodule of the right
$C$-comodule $E_x$, thus $E_x$ is the injective envelope of $S_x$
in the category of right $C$-comodules.
\end{pr}
{\bf Proof:} Let $z=\sum_{x\leq y}\alpha_ye_{x,y}\in E_x$, $z\neq
0$. Choose some $y_0$ such that $\alpha_{y_0}\neq 0$. Then
$p_{x,y_0}\rightarrow z=\alpha_{y_0}e_{x,x}\in S_x\setminus \{
0\}$. \qed

\begin{co}
Any indecomposable injective right $C$-comodule is isomorphic to
$E_x$ for some $x\in X$.
\end{co}

It is easy to describe the left $C^*$-submodule generated by some
$e_{u,v}$.

\begin{pr}
Let $u\leq v$ in $X$. Then $C^*\rightarrow e_{u,v}=< e_{u,z}\; |\;
u\leq z\leq v >$.
\end{pr}
{\bf Proof:} The inclusion $\subseteq$ is clear since
$c^*\rightarrow e_{u,v}=\sum_{u\leq z\leq v}c^*(e_{z,v})e_{u,z}$.
On the other hand for any $u\leq z\leq v$ one has
$p_{z,v}\rightarrow e_{u,v}=e_{u,z}$, so each such $e_{u,z}$ lies
in $C^*\rightarrow e_{u,v}$, and the other inclusion is clear.
\qed

In a similar way, the subspace $S'_x=\sum_{y\in {\cal
C}_i}ke_{y,x}$ is a simple left $C$-subcomodule and
 $E'_x=<e_{y,x}\; |\;
y\in X, y\leq x >$ is its injective envelope as a left
$C$-comodule. Any indecomposable injective left $C$-comodule is
isomorphic to some $E'_x$. There is a good left-right connection
between simple comodules, the next result shows. We note that the
left $C^*$-module structure of $S_v$ induces a right $C^*$-module
structure on $S_v^*$; we denote this right action by $\leftarrow$.
The usual right $C^*$-action on $C$ is denoted by $\leftarrow$,
too.

\begin{pr}
If $v\in X$, then $S^*_v\simeq S'_v$ as right $C^*$-modules.
\end{pr}
{\bf Proof:} We denote by $(\tilde{e}_{v,z})_{z\sim v}$ the basis
of $S_v^*$ dual to $(e_{v,z})_{z\sim v}$. Let
$\gamma:S_v^*\rightarrow S'_v$ be the linear isomorphism such that
$\gamma (\tilde{e}_{v,z})=e_{z,v}$ for any $z\sim v$. We show that
$\gamma $ is a morphism of right $C^*$-modules. Indeed, if $c^*\in
C^*$, then \bea \tilde{e}_{v,z}\leftarrow c^*&=&\sum_{y\sim v}
(\tilde{e}_{v,z}\leftarrow c^*)(e_{v,y})\tilde{e}_{v,y}\\
&=&\sum_{y\sim v} \tilde{e}_{v,z}(c^*\rightarrow
e_{v,y})\tilde{e}_{v,y}\\
&=&\sum_{y\sim v}\sum_{u\sim
v}c^*(e_{u,y})\tilde{e}_{v,z}(e_{v,u})\tilde{e}_{v,y}\\
&=&\sum_{y\sim v}c^*(e_{z,y})\tilde{e}_{v,y} \eea so $\gamma
(\tilde{e}_{v,z}\leftarrow c^*)=\sum_{y\sim
v}c^*(e_{z,y})e_{y,v}=e_{z,v}\leftarrow c^*=\gamma
(\tilde{e}_{v,z})\leftarrow c^*$. \qed

\section{The main result}

We recall that a coalgebra $C$ is called right quasi-co-Frobenius
if $C$ can be embedded as a right $C^*$-module in a free right
$C^*$-module. It is clear that any right co-Frobenius coalgebra is
right quasi-co-Frobenius. It is known that a cosemisimple
coalgebra is right co-Frobenius.

If $(X,\leq )$ is a locally finite preordered set, then the factor
set $\tilde{X}=X/\sim$ with respect to the equivalence relation
$\sim$, associated with $\leq$, is a locally finite partially
ordered set (and there is no danger of confusion if we also denote
the induced order relation on $\tilde{X}$ by $\leq$), so we can
consider the corresponding incidence coalgebra $IC(\tilde{X})$.
Now we can prove our main result. In the particular case where $X$
is a partially ordered set (so then $\tilde{X}=X$), we obtain
\cite[Corollary 3.2]{dnv}. We note that even in this particular
case, our approach is different from the one in \cite{dnv}.

\begin{te} \label{teorema}
Let $C=IC(X)$ be the incidence coalgebra of a locally finite
preordered set $X$. The following assertions are
equivalent.\\
(1) $C$ is right co-Frobenius.\\
(2) $C$ is right quasi-co-Frobenius.\\
(3) $C$ is cosemisimple.\\
(4) For any $x,y\in X$ such that $x\leq y$, we also have that
$y\leq x$ (or equivalently, elements in different equivalence
classes with respect to $\sim$ are not comparable).\\
(5) $IC(\tilde{X})$ is right co-Frobenius.\\
(6) The order relation on $\tilde{X}$ is the equality.
\end{te}
{\bf Proof:} $(1)\Rightarrow (2)$ is clear.\\
$(2)\Rightarrow (3)$ Any indecomposable injective right
$C$-comodule is isomorphic to $E_x$ for some $x\in X$. Take one
such $E_x$. Since $C$ is right quasi-co-Frobenius, $E_x$ is
isomorphic to the dual of  an indecomposable injective left
$C$-comodule $P$ (see \cite[Proposition 1.4]{iovint} or \cite[Proposition 1.3]{iovfrob}),
so it is local (see \cite[Lemma 1.4]{iovcofrobenius}), and then by
Proposition \ref{common} it is generated as a left $C^*$-module by
an element $e_{u,v}$. Thus $E_x=C^*\rightarrow e_{u,v}$.

Then $e_{u,v}\in E_x$, so necessarily $u=x$. Indeed, if $u\neq x$,
then $e_{u,v}$ cannot lie in the span of all $e_{x,y}$, with
$x\leq y$. Thus we have $E_x=<e_{x,z}\;|\; x\leq z\leq v >$,  so
$\{ y\in X\; |\; x\leq y\}=\{ z\;|\; x\leq z\leq v\}$.

Let $\phi:E_x\rightarrow S_v$ be the linear map such that $\phi
(e_{x,s})=e_{v,s}$ for any $s\sim v$, and $\phi (e_{x,s})=0$ for
any $x\leq s<v$. It is clear that $\phi$ is surjective. We show
that $\phi$ is a morphism of left $C^*$-modules. Indeed, if $s\sim
v$ and $c^*\in C^*$, then \bea \phi (c^*\rightarrow
e_{x,s})&=&\sum_{x\leq z\leq
s}c^*(e_{z,s}\phi (e_{x,z}))\\
&=&\sum_{z\sim v}c^*(e_{z,s})e_{v,z}\\
&=&c^*\rightarrow e_{v,s}\\
&=&c^*\rightarrow \phi (e_{x,s}),\eea while if $x\leq s<v$, then
$\phi (c^*\rightarrow e_{x,s})=0$ since $z<v$ for any $x\leq z\leq
s$, and $c^*\rightarrow \phi (e_{x,s})=c^*\rightarrow 0=0$.

In a similar way to what we have done with right $C$-comodules,
when we work with left $C$-comodules, an indecomposable injective
left $C$-comodule is isomorphic to a comodule of the form $E'_w=<
e_{y,w}\; |\; y\in X, y\leq w >$, the indecomposable injective
left $C$-subcomodule of $C$ which contains $e_{w,w}$. Since $C$ is
right quasi-co-Frobenius, it must be right semiperfect (see
\cite[Corollary 3.3.6]{dnr}), so $E'_w$ is finite dimensional (see
\cite[Theorem 3.2.3]{dnr}). Then we have a surjective morphism of
left $C^*$-modules $\theta:(E'_w)^*\rightarrow S_v$, and this
induces an injective morphism $\psi:S_v^*\rightarrow E'_w$ of
right $C^*$-modules. Since the socle of $E'_w$ is just $S'_w$, we
must have $w=v$, so then $E_x\simeq (E'_v)^*$.

Thus we have showed that if $E_x=C^*\rightarrow e_{u,v}$, then
$u=x$, $\{ y\; |\; x\leq y\}=\{ z\; |\; x\leq z\leq v \}$,  and
$E_x\simeq (E'_v)^*$.

Now since $\{ y\; |\; x\leq y\}=\{ z\; |\; x\leq z\leq v \}$, we
see that if $v\leq z$, then $z\sim v$. It follows that
$E_v=S_v=C^*\rightarrow e_{v,v}$. By the previous considerations
for $v$ instead of $x$, we obtain that $E_v\simeq (E'_v)^*$. It
follows that $E_x\simeq E_v$, so then $v=x$. We conclude that
$E_x=S_x$, and then $C=\oplus _{x\in X}E_x=\oplus_{x\in X}S_x$ is
cosemisimple. \\
$(3)\Rightarrow (1)$ holds for any coalgebra.\\
$(3)\Leftrightarrow (4)$ follows from the fact that $S_x$ is the
socle of $E_x$, so $C$ is cosemisimple if and only if $E_x=S_x$
for any $x\in X$.\\
$(4)\Leftrightarrow (6)$ is clear.\\
$(5)\Leftrightarrow (6)$ follows from $(1)\Leftrightarrow (4)$
applied to $(\tilde{X},\leq )$.\qed

Since the cosemisimple property on coalgebras is left-right
symmetric (or alternatively since the conditions (4) and (6) do
not depend on the side we work on), the assertions of Theorem
\ref{teorema} are also equivalent to the left hand versions of
(1), (2) and (5).

Let $n$ be a positive integer and denote by $e_{i,j}$ the matrix
units in the matrix algebra $M_n(k)$. Let ${\mathcal B} \subseteq
\{ 1,\ldots ,n\}\times \{ 1,\ldots ,n\}$ be such that $(i,i)\in
{\mathcal B}$ for any $i\in \{ 1,\ldots ,n\}$, and $(i,k)\in
{\mathcal B}$ whenever $(i,j)\in {\mathcal B}$ and $(j,k)\in
{\mathcal B}$. Thus $\mathcal B$ is a preorder relation on $\{
1,\ldots ,n\}$. Then ${\mathcal M}({\mathcal B},k)=\sum_{(i,j)\in
{\mathcal B}}ke_{i,j}$ is a subalgebra of ${\mathcal M}_n(k)$,
called the structural matrix algebra associated to $\mathcal B$.
It  consists of all matrices whose $(i,j)$-entries are zero for
$(i,j)\notin {\mathcal B}$. The dual coalgebra of ${\mathcal
M}({\mathcal B},k)$ is just the incidence coalgebra of the
preordered set $\{ 1,\ldots ,n\}$, with the relation defined by
$\mathcal B$. Since a finite dimensional algebra is Frobenius if
and only if its dual coalgebra is right co-Frobenius (and then by
the finite dimensionality it is also left co-Frobenius), we obtain
as a consequence of Theorem \ref{teorema} the following.

\begin{co} \label{corollaryFrobeniusstructural}
A structural matrix algebra ${\mathcal M}({\mathcal B},k)$ is
Frobenius if and only if ${\mathcal B}=(I_1\times I_1)\cup \ldots
\cup (I_r\times I_r)$ for some partition $I_1,\ldots ,I_r$ of $\{
1,\ldots ,n\}$. In this case ${\mathcal M}({\mathcal B},k)\simeq
M_{n_1}(k)\times \ldots \times M_{n_r}(k)$, where
$n_1=|I_1|,\ldots,n_r=|I_r|$.
\end{co}

\begin{re}
It is possible to prove Corollary
\ref{corollaryFrobeniusstructural} directly, using only ring
theory methods. The referee indicated us one such proof. A
structural matrix algebra may be put (up to an isomorphism) into
upper block triangular form, see for instance \cite[page 28]{abw}.
As a consequence of \cite[Chapter III, Proposition 2.7]{ars}, the
global dimension of a structural matrix algebra is finite. On the
other hand, the global dimension of a Frobenius algebra is either
0 or $\infty$ by \cite[Theorem 11]{en} or \cite[Proposition
15]{a}. Then if a structural matrix algebra is Frobenius, it has
global dimension 0, so it must be semisimple, and then it consists
only of diagonal blocks.
\end{re}

We end by presenting a categorical connection between the two
incidence coalgebras, $C=IC(X)$ and $IC(\tilde{X})$. More
precisely, we show that these two coalgebras are Morita-Takeuchi
equivalent, i.e. their categories of right comodules are
equivalent. Let $\cal S$ be a system of representatives for the
equivalence classes of $X$ with respect to $\sim$, and let $m\in
C^*$ be such that $m(e_{u,v})=1$ if $u=v\in {\cal S}$, and
$m(e_{u,v})=0$ for any other $u\leq v$. Thus $m$ is just $\eps$ on
$E_u$ with $u\in{\cal S}$, and $m$ is zero on any other $E_u$.
Then $m$ is an idempotent of $C^*$, and by \cite[Lemma 6]{rad} we
have that $m\rightarrow C\leftarrow m$ has a coalgebra structure
with the comultiplication defined by
$$\Delta' (m\rightarrow c\leftarrow m)=\sum (m\rightarrow c_1\leftarrow
m)\otimes (m\rightarrow c_2\leftarrow m)$$ and counit just the
restriction of $\eps$ to $m\rightarrow C\leftarrow m$. Note that
$m\rightarrow C\leftarrow m$ is not a subcoalgebra, but a factor
coalgebra of $C$. Using the concept of a basic coalgebra defined
in \cite{cm}, and the terminology and results of \cite[Section
3]{cg}, we have that $m$ is a basic idempotent of $C$ and
$m\rightarrow C\leftarrow m$ is the basic coalgebra of $C$, so
then $m\rightarrow C\leftarrow m$ is Morita-Takeuchi equivalent to
$C$ see \cite[Corollary 2.2]{cm} or \cite[Proposition 3.6]{cg}.

By the way $m$ is defined, it is easy to see that $m\rightarrow
e_{u,v}\leftarrow m=e_{u,v}$ if $u,v\in {\cal S}$ and $u\leq v$,
and $m\rightarrow e_{u,v}\leftarrow m=0$ in any other case. Thus
$m\rightarrow C\leftarrow m=< e_{u,v}\;| \; u,v\in {\cal S}, u\leq
v>$, and its comultiplication works as \bea \Delta'
(e_{u,v})&=&\sum_{u\leq y\leq v} (m\rightarrow e_{u,y}\leftarrow
m)\otimes (m\rightarrow e_{y,v}\leftarrow m)\\
&=&\sum_{u\leq y\leq v, y\in {\cal S}} e_{u,y}\otimes e_{y,v}\eea
This shows that the linear map $f:(m\rightarrow C\leftarrow
m)\rightarrow IC(\tilde{X})$ defined by
$f(e_{u,v})=e_{\overline{u},\overline{v}}$ for any $u,v\in {\cal
S}$ with $u\leq v$, is an isomorphism of coalgebras. Here we
denoted by $\overline{u}$ the class of $u$ in the factor set
$X/\sim$. We conclude that $IC(X)$ and $IC(\tilde{X})$ are
Morita-Tacheuchi equivalent.

As a consequence, we obtain by duality that for a finite
preordered set $X$, the incidence algebra of $X$ is Morita
equivalent to the incidence algebra of the partially ordered set
$\tilde{X}$. Thus a structural matrix algebra, which is isomorphic
to a blocked matrix algebra via a permutation of rows and columns,
is Morita equivalent to a more simple structural matrix algebra,
which is presented in the block form with blocks of size 1.

{\bf Acknowledgment.} We would like to thank the referee for the
valuable remarks and suggestions. The research of the first two
authors was supported by the UEFISCDI Grant
PN-II-ID-PCE-2011-3-0635, contract no. 253/5.10.2011 of CNCSIS.
The third author was supported by the Sectorial Operational
Programme Human Resources Development (SOP HRD), financed from the
European Social Fund and by the Romanian Government under the
contract number SOP HRD/107/1.5/S/82514.

\end{document}